\documentclass[draft]{amsart}
\usepackage{latexsym}
\usepackage{amsfonts}
\usepackage{amsmath}
\usepackage{fancyheadings}
\usepackage{amssymb}
\usepackage{graphics}
\usepackage[OT2,OT1]{fontenc}
\usepackage{wrapfig}
\usepackage{color}
\newtheorem{theorem}{Theorem}[section]

\newtheorem{lemma}[theorem]{Lemma}

\theoremstyle{definition}

\newcommand{\der}[2]{#1^{(#2)}}

\newcommand{\lcs}[2]{\gamma_{#2}(#1)}

\newcommand{\sym}{\sf S}

\newcommand{\ord}[1]{|\,#1\,|}

\newcommand\rc[3]{[#1,\,{}^{#3}\!#2]}

\renewcommand{\d}{d}

\renewcommand{\leq}{\leqslant}
\renewcommand{\geq}{\geqslant}

\begin{document}

\title{Small derived quotients in finite $p$-groups}
\author{Csaba Schneider}
\address{Informatics Laboratory\\ 
Computer and Automation
Research Institute\\
The Hungarian Academy of Sciences\\
1518 Budapest Pf.\ 63\\
Hungary}
\email{csaba.schneider@sztaki.hu\protect{\newline} {\it WWW:}
www.sztaki.hu/$\sim$schneider}

\begin{abstract}
More than 70 years ago, P.\ Hall showed that if $G$ is a finite $p$-group
such that a term $\der G{d+1}$ of the
derived series is non-trivial, then the order of the quotient
$\der Gd/\der G{d+1}$ is
at least
$p^{2^d+1}$.
Recently Mann proved that, in a finite $p$-group,
Hall's lower bound can be taken for at most two distinct $d$. 
For odd $p$, we prove a sharp version of 
this result and characterise the groups with two small derived quotients.
\end{abstract}

\date{\today}
\maketitle
\dedicatory{Dedicated to the memory of my dear friend and mentor, Edit Szab\'o.}

\section{Introduction}
Suppose that $G$ is a finite $p$-group in which a term
$\der G{d+1}$ of the derived series is non-trivial 
(we index the terms of the derived series so that
$\der G0=G$, $\der G1=G'$, $\der G2=G''$, etc). Then how small can  
the order of the
quotient $\der Gd/\der G{d+1}$ possibly be? 
As far as I know, the answer for this general question
is not known. 
Hall showed in~\cite{Hall34} that if  
$H$ is a non-abelian normal subgroup in a finite $p$-group $G$
that is contained in the $i$-th term $\gamma_i(G)$ of the lower central series
of $G$ (the terms of the lower central series are indexed so that $\lcs G1
=G$, $\lcs G2=G'$, etc), then $|H/H'|\geq p^{i+1}$ (see Lemma~\ref{lemma}(a)).
As $\der Gd\leq\lcs G{2^d}$,
this result implies that $\log_p|\der Gd/\der G{d+1
}|
\geq {2^d+1}$ provided $\der G{d+1}\neq 1$.

In a finite $p$-group $G$ let us call a quotient $\der Gd/\der G{d+1}$ a {\em
  small derived quotient} if $\der G{d+1}\neq 1$ and $\log_p
|\der Gd/\der G{d+1}|=2^d+1$. Mann~\cite{Mann} showed that a finite $p$-group
  can have at most two small derived quotients. Building on the results
  of~\cite{sch2}, we prove the following 
sharp theorem.

\begin{theorem}\label{main}
Let $p$ be an odd prime and let $G$ be a finite $p$-group that contains two
small derived quotients. Then $p\geq 5$, $|G|=p^6$, $|G''|=p$, and $G$ has
class~$5$. Further, for $p\geq 5$, there are precisely
$p+4+\gcd(4,p-1)+\gcd(5,p-1)+\gcd(6,p-1)$ pairwise non-isomorphic finite $p$-groups with
two small derived quotients.

\end{theorem}

My main motivation for studying small derived quotients in 
$p$-groups was to improve the existing lower bounds for the order of a
$p$-group with a given derived length $d+1$. In such a group $\der G{d}\neq
1$. If we assume, as did Philip Hall in~\cite{Hall34}, that, for 
$i=0,\ldots,d-1$, the quotient $\der Gi/\der G{i+1}$ is small, then we obtain
that $\log_p|G|\geq 2^d+d$. However, if we use Mann's result that at most two
of the derived quotients can be small, we find $\log_p|G|\geq
2^d+2d-2$; see~\cite{Mann}. 
Using Theorem~\ref{main} we can easily obtain a miniscule
improvement of Mann's lower bound for $|G|$. However, in a separate
article~\cite{sch2}, I show that investigating the metabelian quotients of
$G$, the linear term in Mann's bound can be further improved. To be precise, 
if $p\geq 5$ and 
$\der G{d}\neq 1$, then $|G|\geq 2^d+3d-6$; see~\cite{sch2} for
details.

\section{The structure of small derived quotients}

If $A$ and $B$ are subgroups in a group $G$ 
and $n$ is a natural
number then let $\rc{A}{B}{n}$ denote the left-normed commutator subgroup
$$
\rc{A}{B}{n}=[A,\underbrace{B,\ldots,B}_{n\rm\ copies}].
$$
One can easily show by induction on $i$ that 
if $A$ and $B$ are normal subgroups of $G$,  then
\begin{equation}\label{repcomm}
[A,\lcs Bi]\leq \rc ABi.
\end{equation}

We will need the following well-known lemma. Part~(a) was shown in
~\cite{Hall34}, while part~(b) can be found as~\cite[Lemma~2.1]{max-class}.

\begin{lemma}\label{lemma}
(a) Suppose that $H$ is a non-abelian normal subgroup in a finite $p$-group $G$
such that $H\leq \lcs Gi$. Then $|H/H'|\geq p^{i+1}$ and $|H|\geq p^{i+2}$. 

(b) If $G$ is a group and $H$ is a normal subgroup such that $G/H$ is cyclic,
then $G'= [G,H]$. 
\end{lemma}

Suppose that $G$ is a finite $p$-group and that 
$\der Gd/\der G{d+1}$ is a small derived quotient for some $d\geq 0$. 
As $\der Gd\leq\lcs G{2^d}$, we obtain
$$
\der G{d+1}
=[\der Gd,\der Gd]\leq [
\der Gd,\lcs{G}{2^\d}]\leq \rc{\der Gd}{G}{2^\d},
$$
therefore  we have
the following chain of $G$-normal subgroups:
\begin{equation}\label{chain}
\der Gd>[\der Gd,G]>[\der Gd,G,G]>\cdots>\rc {\der Gd}{G}{2^\d}\geq \der G{d+1}.
\end{equation}
Counting number of non-trivial factors of this chain, we obtain that $\der Gd/[\der Gd,G]$ has order at most $p^2$.
If $\der Gd/[\der Gd,G]$ is cyclic, then, by Lemma~\ref{lemma}(b),
the subgroup $\der G{d+1}$ coincides with $[\der Gd,[\der Gd,G]]$, and so 
$$
\der G{d+1}
=[\der Gd,[\der Gd,G]]\leq [
\der Gd,\lcs{G}{2^\d+1}]\leq \rc{\der Gd}{G}{2^\d+1}.
$$
Thus, in this case, we
obtain the following modified chain:
\begin{equation}\label{newchain}
\der Gd>[\der Gd,G]>[\der Gd,G,G]>\cdots>\rc{\der Gd}{G}{2
^d+1}\geq\der G{d+1}.
\end{equation}
If the first quotient in these chains has order $p$, then this quotient is
cyclic, and so~\eqref{newchain} must hold. 
In this case, counting the non-trivial factors in~\eqref{newchain}, we find 
that the following chain must be valid:
\begin{equation}\label{ch1}
\der Gd>[\der Gd,G]>[\der Gd,G,G]>\cdots>
\rc{\der Gd}{G}{2
^d+1}=\der G{d+1}.
\end{equation}
Now suppose that the 
first quotient $\der Gd/[\der Gd,G]$ has order $p^2$. 
Then chain~\eqref{newchain} is too long, and
so $\der Gd/[\der Gd,G]$ must be elementary abelian. 
As before, we count the number of factors
in~\eqref{chain} and find the following chain:
\begin{equation}\label{ch2}
\der Gd>[\der Gd,G]>[\der Gd,G,G]>\cdots>\rc{\der Gd}{G}{2^\d}=\der G{d+1}.
\end{equation}

It is, perhaps, somewhat surprising that, in general, chain~\eqref{ch2} is not
possible. 

\begin{theorem}\label{smallth}
Suppose that $p$ is an odd prime,  $d\geq 1$, and that $\der Gd/\der G{d+1}$ is a 
small derived quotient in a finite $p$-group $G$. Then 
$|\der Gd/[\der Gd,G]|=p$ and so chain~$\eqref{ch1}$ must be valid.
\end{theorem}

Theorem~\ref{smallth} first appeared in my PhD thesis~\cite{csthesis}. The
special case of  $d=1$ was also proved in a recent article~\cite{sch}. The proof
of the general case can be found, besides my thesis, in the forthcoming
article~\cite{sch2}.  

\section{Proof of Theorem~\ref{main}}

Let $p$ be an odd prime,
let $G$ be a finite $p$-group and let $d$ be a non-negative 
integer such that 
$\der Gd/\der G{d+1}$ is a small derived quotient. Let us assume, in addition, 
that
$d$ is the smallest such integer. If~\eqref{ch1} is valid, 
then 
$$
\der G{d+1}\leq \rc{\der Gd}G{2^d+1}\leq \lcs G{2^{d+1}+1}.
$$
Now easy induction shows, for $e\geq 1$, that $\der G{d+e}\leq\lcs
G{2^{d+e}+2^{e-1}}$. Hence Lemma~\ref{lemma}(a) implies that $\der G{d+e}/\der
G{d+e+1}$ cannot be small for $e\geq 1$. Therefore, in this case, $\der
Gd/\der G{d+1}$ is the unique small derived quotient in $G$.

Suppose now that~\eqref{ch2} is valid. In this case, it is easy to show that
$\der G{d+1}/[\der G{d+1},G]$ must be cyclic (see~\cite[Corollary~5.2]{sch2}), 
and following the argument in
the previous paragraph, one easily obtains that $\der G{d+e}/\der G{d+e+1}$
cannot be small for $e\geq 2$. Hence only the 
derived quotients $\der
Gd/\der G{d+1}$ and $\der G{d+1}/\der G{d+2}$ can be small in $G$. By assumption, 
\eqref{ch2} must hold
for the quotient 
$\der Gd/\der G{d+1}$ and, as shown above, \eqref{ch1}~must be valid for the quotient $\der G{d+1}/\der
G{d+2}$. 

So far, we have obtained Mann's result in~\cite{Mann} that a finite
$p$-group can have at most two small derived quotients (the assumption that
$p$ is odd has played no r\^ole up to this point). Now we may use
Theorem~\ref{smallth} and obtain, for $p\geq 3$, that~\eqref{ch2} is only
possible for $d=0$. Thus if $G$ has odd order, then
the two distinct small derived quotients must be $G/G'$, $G'/G''$.
The quotient $G'/G''$ is as in~\eqref{ch1} and so we find that
$G''=[G',G,G,G]=\lcs G5$.
As 
$|G'/G''|=p^3$, a result that Blackburn attributes 
to P.\ Hall (see~\cite{black87})
shows that $|G''|=p$. Thus $|G|=p^6$, and, as $G''=\lcs G5\neq 1$, we obtain
that $G$ has class~5. Therefore $G$ is a group with maximal class.

It remains to show that 
the restriction on $p$ in the theorem holds and that
the number of groups with two small derived quotients is as claimed. We still
work under the assertion that $p$ is odd and that 
$G$ has two small derived quotients.  As chain~\eqref{ch1} is valid for $G'/G''$, we obtain
$$
[\lcs G2,\lcs G3]=[G',[G',G]]=G''=\lcs G5,
$$
and so $G$ has degree of commutativity~0 (see~\cite[page~57]{max-class}). A 3-group with two distinct small
derived quotient lies in Blackburn's class ${\sf ECF}(6,6,3)$ and 
so~\cite[Theorem~3.8]{max-class} shows that such a $3$-group has degree of
commutativity greater than zero. Thus we obtain that $p\geq 5$. (The claim that
$p\geq 5$ can also be verified using the Small Groups Library of the
computational algebra systems~\cite{GAP} or~\cite{Magma}.)

Let $H$ be a $p$-group of maximal class with order $p^6$. As $H'/[H',H]$ is
cyclic with order $p$, we obtain that $H''\leq\lcs H5$ (Lemma~\ref{lemma}(b)). 
Thus, by the above, 
$H$ has two distinct small derived quotients, if and only if $H$ is not
metabelian. By~\cite[Theorems~4.4 and 4.5]{max-class}, the number of such
groups is $p+4+\gcd(4,p-1)+\gcd(5,p-1)+\gcd(6,p-1)$.

Thus the proof of Theorem~\ref{main} is now complete.

\section{Some final remarks}

The Sylow
2-subgroup $P$ of the symmetric group $\sym_{2^\d}$ of rank $2^\d$ satisfies
$$
\log_2\ord{\der{P}{\d-2}/\der{P}{\d-1}}=2^{\d-2}+1
\quad\mbox{and}\quad \der{P}{\d-1}\neq 1
$$
(see \cite[Lemma~(II.7)]{lg}). 
Hence the derived quotient
$\der{P}{\d-2}/\der{P}{\d-1}$ is small, and one can also show
using~\cite[Lemma~(II.7)]{lg} that, in this case, ~\eqref{ch2} is valid; that
is, $|\der{P}{d-2}/[\der{P}{d-2},P]|=p^2$. Therefore Theorem~\ref{smallth} is not valid for $2$-groups

There are many finite $p$-groups in which the quotient $G/G'$ is small. 
Finite $p$-groups in which $G'/G''$ is small were
characterised 
in~\cite{sch}. However, for odd $p$, it is not clear whether
in a $p$-group $G$ the quotient $\der Gd/\der G{d+1}$ can be small for $d\geq
2$. We do not even know of odd-order examples $G$ in
which $\der G2/\der G3$ is small, that is, $\der G3\neq 1$ and $|\der G2/\der
G3|=p^5$.

\section*{Acknowledgment}
Much of the research presented in this paper was carried out while I was a PhD
student at The Australian National University. I am particularly grateful to
my PhD supervisor, Mike Newman, for his continuous support.

\bibliographystyle{alpha}

\end{document}